\documentclass{amsart}
\usepackage{amsfonts}
\usepackage{amsmath,amssymb}
\usepackage{amsthm}
\usepackage{amscd}
\usepackage{graphics}
\usepackage{graphicx}

\theoremstyle{definition}{
\newtheorem{Def}{{\rm Definition}}
\newtheorem{Ex}{{\rm Example}}
\newtheorem{Rem}{{\rm Remark}}

}
\theoremstyle{plain}
{

\newtheorem{Prop}{Proposition}
\newtheorem{Thm}{Theorem}
\newtheorem{MainThm}{Main Theorem}

}

\begin{document}
\title[New round fold maps on some $7$-dimensional manifolds into the plane]{Round fold maps into the plane on some $7$-dimensional closed and simply-connected manifolds}
\author{Naoki Kitazawa}
\address{Institute of Mathematics for Industry, Kyushu University, 744 Motooka, Nishi-ku Fukuoka 819-0395, Japan\\
 TEL (Office): +81-92-802-4402 \\
 FAX (Office): +81-92-802-4405 \\
}
\email{n-kitazawa@imi.kyushu-u.ac.jp}
\urladdr{https://naokikitazawa.github.io/NaokiKitazawa.html}
\maketitle
\begin{abstract}
{\it Round} fold maps are smooth maps on closed manifolds which are locally represented as the product maps of Morse functions and identity maps on open disks and whose singularity is realized as concentrically embedded spheres. The author previously introduced such maps. Our paper presents round fold maps on some $7$-dimensional simply-connected manifolds whose cohomology rings are isomorphic to that of the product of the $2$-dimensional complex projective space and a $3$-dimensional sphere. Such manifolds have been studied precisely by Wang and round fold maps on spin manifolds in these manifolds have been previously studied by the author. These manifolds form explicit classes of higehr dimensional closed and simply-connected manifolds, which are central objects in classical algeberic topology and differential topology. Understanding these manifolds in geometric and constructive ways is still attractive, which we think as pioneers.

{\it Fold} maps are defined as smooth maps which are locally represented as the product maps of Morse functions and identity maps on open disks. They are fundamental and strong tools in generalizations of theory of Morse functions and applications to geometry of manifolds. Explicit construction of fold maps are difficult even on elementary or well-known manifolds whereas we can know the (non-)existence from Eliashberg's celebrating theory in the 1970s and related one in considerable cases.

\end{abstract}

% Leave these items like this, and the journal will fill them in.
% \received{Month Day, Year}   % receive date (for example: October 11, 1999)
% \revised{Month Day, Year}    % date of revision; omit, if no revision;
%                             % if multiple revisions, separate by commas
% \published{Month Day, Year}  % publish date\submitted{Bill Murray}      % Name of Journal's Editor, 
% who handled Article 
% \volumeyear{2014} % Volume Year
% \volumenumber{16} % Volume Number 
% \issuenumber{2}   % Issue Number
% \startpage{1}     % PageNumber of first page
% \articlenumber{1} % Sequence number of article within issue
% If copyright is retained by author, comment this out:
% \owner{International Press}

\maketitle
\section{Introduction, terminologies and notation.}
\label{sec:1}

{\it Fold} maps form an important class of smooth maps containing the class of Morse functions. They play important roles in singularity theory of smooth maps and applications to geometry, especially, algebraic topology and differential topology of manifolds. \cite{golubitskyguillemin} is for systematic expositions on Morse functions, fold maps, and more general smooth maps from the viewpoint of singularity theory of smooth maps. For fold maps and applications to algebraic topological and differential topological understandings of manifolds, we review when we need. 

For Morse theory or theory of Morse funtions, see \cite{milnor2} and \cite{milnor3}. \cite{milnor2} explains about fundamental philosophy and applications to investigating the spaces of important paths such as geodesics on Riemannian manifolds and Lie groups and the spaces of such paths on them. \cite{milnor3} explains about elementary tools in applying our theory of Morse functions to investigate algebraic topological and differential topological properties of the manifolds. This also presents fundamental strong facts such as the correspondence between so-called {\it handles} and singularities of Morse functions and such facts are important in the present paper. We omit rigorous expositions on them.  

Here we introduce some fundamental terminologies and notation.

For an integer $k>0$, ${\mathbb{R}}^k$ denotes the $k$-dimensional Euclidean space and we regard this as a natural smooth manifold and a Riemannian manifold endowed with the standard Euclidean metric. $||x|| \geq 0$ denotes the distance between $x \in {\mathbb{R}}^k$ and the origin $0 \in {\mathbb{R}}^k$. $S^k:=\{x \in {\mathbb{R}}^{k+1} \mid ||x||=1\}$ denotes the $k$-dimensional unit sphere for an integer $k \geq 0$ and $D^k:=\{x \in {\mathbb{R}}^{k} \mid ||x|| \leq 1\}$ denotes the $k$-dimensional unit disk for an integer $k \geq 0$. They are $k$-dimensional smooth submanifolds in ${\mathbb{R}}^k$. The unit sphere is a closed submanifold with no boundary and the unit disk is a closed submanifold whose boundary is not empty.

The $1$-dimensional Euclidean space ${\mathbb{R}}^1$ is usually denoted by $\mathbb{R}$. $\mathbb{Z} \subset \mathbb{R}$ denotes the ring of all integers and $\mathbb{N} \subset \mathbb{Z}$ denotes the set of all positive integers or integers greater than $0 \in \mathbb{Z}$.

A {\it homotopy sphere} means a smooth manifold homeomorphic to a unit sphere. It is {\it standard} ({\it exotic}) if it is diffeomorphic to a unit sphere (resp. not diffeomorphic to any unit sphere).

A {\it singular point} $p$ of a smooth map $f:M \rightarrow N$ means a point in the manifold of the domain where the rank of the differential ${df}_p$ of the manifold is both the dimensions of the manifolds $M$ and $N$.
The {\it singular set} $S(f)$ is the set of all singular points of $f$. A {\it singular value} of $f$ is a point in the image $f(S(f))$, called the {\it singular value set} of $f$.

\begin{Def}
	\label{def:1}
Let $m \geq n \geq 1$ be integers.
A smooth map $f:M \rightarrow N$ from an $m$-dimensional manifold $M$ with no boundary into an $n$-dimensional manifold $N$ with no boundary is said to be a {\it fold} map if on some open small neighborhood of each singular point $p$, $f$ has the form
$$(x_1,\cdots,x_m) \rightarrow (x_1,\cdots,x_{n-1},{\Sigma}_{j=1}^{m-n-i(p)+1} {x_{n+j-1}}^2-{\Sigma}_{j=1}^{i(p)} {x_{m-i(p)+j}}^2)$$ for suitable local coordinates and a suitable integer $0 \leq i(p) \leq \frac{m-n+1}{2}$. 
\end{Def}
\begin{Prop}
	In Definition \ref{def:1}, the singular set $S(f)$ is a smooth closed submanifold of dimension $n-1$ with no boundary. Furthermore, the restriction $f {\mid}_{S(f)}$ is a smooth immersion.
	$i(p)$ is unique and it is called the {\rm index} of $p${\rm :} note also that it is also called the {\rm absolute index} of $p$ however we adopt "the index".
\end{Prop}

For example, as an important fact, the existence of fold maps is an important study, discussed first in \cite{thom,whitney} and later \cite{levine}. A closed manifold admits a fold map into ${\mathbb{R}}^2$ if and only if the Euler number is even according to these studies. Remember that Morse functions always exist plentifully. \cite{eliashberg,eliashberg2} are celebrating studies on the existence of fold maps by Eliashberg. It is also important to know that explicit construction of Morse functions and fold maps is very difficult in general although we know the existence. It lies on the main theme of our paper.

Investigating restrictions on manifolds admitting fold maps of certain good classes, characterizations of certain classes of manifolds by the existence of fold maps of certain good classes and construction of explicit fold maps on elementary manifolds have been attractive and challenging problems. They produce important problems and developments both in singularity theory of smooth maps and applications to geometry, especially, algebraic topology and differential topology of manifolds.

A {\it round} fold map is a fold map whose singular value set is the disjoint union of spheres embedded concentrically into $N:={\mathbb{R}}^n$ where $n \geq 2$. A canonical projection of an $m$-dimensional unit sphere to ${\mathbb{R}}^n$ is a round fold map for $m \geq n \geq 2$ for example. This was first introduced by the author in \cite{kitazawa0.1,kitazawa0.2,kitazawa0.3}.

Main Theorem \ref{mthm:1} is our main result.
 This can be regarded as an extension of Theorem \ref{thm:6} in the next section.
 This is shown in \cite{kitazawa3} as a new result by the author.

For a topological space $X$, $H_j(X;\mathbb{Z})$ ($H^j(X;\mathbb{Z})$) denotes the {\it $j$-th} (resp. {\it co}){\it homology group} whose {\it coefficient ring} is $\mathbb{Z}$ or the {\it $j$-th integral} (resp. {\it co}){\it homology group}. The direct sum ${\oplus}_{j=0}^{\infty} H^j(X;\mathbb{Z})$ denotes the direct sum of the $j$-th cohomology groups for all non-negative integers $j$. It is also denoted by $H^{\ast}(X;\mathbb{Z})$. It is regarded as a graded commutative ring where the product is given by the {\it cup product} $u_1 \cup u_2$ for $u_1,u_2 \in H^{\ast}(X;\mathbb{Z})$. This is the {\it cohomology ring} of $X$ whose {\it coefficient ring} is $\mathbb{Z}$ or the {\it integral cohomology ring} of $X$. 
${\mathbb{C}P}^k$ denotes the $k$-dimensional complex projective space, which is a $k$-dimensional complex manifold and $2k$-dimensional closed and simply-connected smooth manifold.
\begin{MainThm}
	\label{mthm:1}
	Let $m>6$ be an integer.
There exists a family $\{M_{j,s}\}_{(j,s) \in \mathbb{Z} \times \{0,1\}}$ of $m$-dimensional closed and simply-connected manifolds whose integral cohomology rings are isomorphic to that of ${\mathbb{C}P}^2 \times S^{m-4}$ and which admit round fold maps into ${\mathbb{R}}^2$ and enjoy the following properties.
\begin{enumerate}
	\item $M_{j,0}$ and $M_{j,1}$ are not homeomorphic for every integer $j$.
	\item $M_{j_1,s_1}$ and $M_{j_2,s_2}$ are not homeomorphic unless $(j_1,s_1)=(\pm j_2,s_2)$. 
\end{enumerate}
\end{MainThm}

 Including these manifolds, closed and simply-connected manifolds have been central objects in algebraic topology and differential topology of manifolds since the midst of the 20th century. Such manifolds whose dimensions are at least $5$ have been classified in algebraic and abstract ways. \cite{ranicki,wall2} show systematic theory on this. $5$, $6$ and $7$-dimensional cases have more explicit classifications via more explicit algebraic topological ways, presented in \cite{barden,crowleyescher,crowleynordstrom,kreck,wang} for example. Related expositions are also presented again in Remark \ref{rem:3}. $4$-dimensional cases are due to Freedman (\cite{freedman}) for example where there exist very difficult unsolved problems on differentiable structures. The $3$-dimensional case is well-known as the problem of Poincar\'e, solved by Perelman. 

In the next section, we review existing studies on fold maps. For example, some examples on manifolds in the following.
\begin{itemize}
	\item The total spaces of bundles over standard spheres whose fibers are closed manifolds and manifolds represented as connected sums of the total spaces of bundles over unit spheres whose fibers are homotopy spheres.
	\item {\it Graph manifolds}, forming an important class of $3$-dimensional closed and orientable manifolds.
	\item $7$-dimensional closed and simply-connected manifolds whose cohomology rings are isomorphic to that of ${\mathbb{C}P}^2 \times S^3$. Theorem \ref{thm:6} presents round fold maps into ${\mathbb{R}}^4$ on these manifolds which are also {\it spin}. We introduce {\it spin} manifolds later shortly and Main Theorem \ref{mthm:1} presents round fold maps into ${\mathbb{R}}^2$ on these manifolds which may not be spin. 
\end{itemize}
 The third section discusses Main Theorem \ref{mthm:1} and related new observations such as Theorems \ref{thm:7} and \ref{thm:8}.
\section{Reviewing studies of fold maps, especially, characterizations of manifolds admitting fold maps of fixed classes.}
 In the present paper, a diffeomorphism between two smooth manifolds means a smooth map which is a homeomorphism and which has no singular points. A diffeomorphism on a fixed smooth manifold is defined in a similar way.

Bundles are fundamental tools here and in geometry. For systematic studies, see \cite{milnorstasheff,steenrod} for example. We assume knowledge on fundamental terminologies, notions and properties on bundles, presented in such books.

A {\it smooth} bundle means a bundle whose fiber is a smooth manifold and whose structure group consists of diffeomorphisms. 

A {\it linear} bundle means a bundle whose fiber is a Euclidean space, a unit sphere, or a unit disk and whose structure group consists of diffeomorphisms regarded as linear transformations in canonical ways. 

%A {\it complex vector bundle} means a bundle satisfying the following conditions.
%%\begin{itemize}
%	\item The fiber is a Euclidean space whose dimension is even and endowed with the standard complex structure, making the fiber the natural complex vector space.
%	\item The structure group consists of linear transformations on the complex vector space where the field is the complex number field.   
%\end{itemize}

%It is also a linear bundle, defined just before.

\begin{Def}
	\label{def:2}
	Let $m \geq n \geq 2$ be integers.
	A fold map $f:M \rightarrow {\mathbb{R}}^n$ on an $m$-dimensional closed and connected manifold $M$ into ${\mathbb{R}}^n$ is said to be {\it round} if there exist a diffeomorphism $\phi:{\mathbb{R}}^n \rightarrow {\mathbb{R}}^n$ and an integer $l>0$ enjoying the relation $\phi(f(S(f)))=\{x \in {\mathbb{R}}^n\mid ||x|| \in \mathbb{N}, ||x|| \leq l\}$. 
\end{Def}
\begin{Ex}
	\label{ex:1}
	The canonical projection of the unit sphere $S^m \subset {\mathbb{R}}^{m+1}$ to ${\mathbb{R}}^n$ mapping $(x_1,x_2) \in {\mathbb{R}}^{m+1}={\mathbb{R}}^n \times {\mathbb{R}}^{m-n+1}$ to $x_1 \in {\mathbb{R}}^{n}$ is a round fold map where we identify ${\mathbb{R}}^{m+1}$ with ${\mathbb{R}}^n \times {\mathbb{R}}^{m-n+1}$ canonically. This is a kind of elementary exercises on the theory of Morse functions and more general theory on singularities of smooth maps.
\end{Ex}
\begin{Rem}
	\label{rem:1}
	For example, \cite{kitazawa0.5} defines a {\it round} fold map in the case ${\mathbb{R}}^n=\mathbb{R}$ as a Morse function obtained by attaching two copies of a Morse function on a compact manifold with no boundary. A Morse function with two singular points on a homotopy sphere, or a Morse function for Reeb's theorem, is round, according to the definition. We review Reeb's theorem in Theorem \ref{thm:3} (\ref{thm:3.1}).
	As a kind of exercises, we can know the following two.
	
\begin{enumerate}
	\item A manifold admitting a round fold map is obtained by attaching two copies of a compact and connected manifold along the boundaries by a suitable diffeomorphism.
	\item If an $m$-dimensional admits a round fold map into ${\mathbb{R}}^n$ where $m \geq n \geq 2$, then it also admits one into ${\mathbb{R}}^{n^{\prime}}$ for any integer $1 \leq n^{\prime} \leq n$. 
	\end{enumerate}
\end{Rem}

	In Definition \ref{def:2}, we can see that the composition of the restriction of $\phi \circ f$ to the preimage ${(\phi \circ f)}^{-1}(\{x \in {\mathbb{R}}^n \mid ||x|| \geq \frac{1}{2}\})$ with the canonical projection onto $S^{n-1}$ mapping $p \in \{x \in {\mathbb{R}}^n \mid ||x|| \geq \frac{1}{2}\}$ to $\frac{1}{||p||}p \in S^{n-1}$ gives a smooth bundle. We call this a {\it bundle for a global monodromy} of $f$. We can also see that the composition of the restriction of $\phi \circ f$ to the preimage ${(\phi \circ f)}^{-1}(\{x \in {\mathbb{R}}^n \mid l^{\prime}-\frac{1}{2} \leq ||x|| \leq l^{\prime}+\frac{1}{2}\})$ with the canonical projection onto $S^{n-1}$ mapping $p \in {(\phi \circ f)}^{-1}(\{x \in {\mathbb{R}}^n \mid l^{\prime}-\frac{1}{2} \leq ||x|| \leq l^{\prime}+\frac{1}{2}\})$ to $\frac{1}{||p||}p \in S^{n-1}$ gives a smooth bundle for each integer $1 \leq l^{\prime} \leq l$. We call such a family of $l$ bundles a {\it family of bundles for componentwise monodromies}.
	 \begin{Def}
	 	Let $f$ be a round fold map. If a bundle for a global monodromy of $f$ is trivial, then $f$ is said to {\it have a globally trivial monodromy}. If there exists a family of bundles for componentwise monodromies of $f$ each of which is a trivial smooth bundle, then $f$ is said to {\it have a family of componentwisely trivial monodromies}.
	 	\end{Def}
 	These notions have been introduced first by the author in \cite{kitazawa0.2,kitazawa0.3,kitazawa0.4,kitazawa0.5} for example where they have been named in different ways.
 	
\begin{Ex}
	\label{ex:2}
	The canonical projection of the unit sphere presented in Example \ref{ex:1} has a globally trivial monodromy.
\end{Ex} 
For a smooth map $f$, the {\it regular value set} of $f$ is the complementary set of the singular value set of $f$. A {\it regular value} of $f$ means a point in the manifold of the target which is not a singular value of $f$. In other words, a regular value means a point in the regular value set. 	
\begin{Thm}[\cite{kitazawa0.4}]
	\label{thm:1}
	Let $m>n \geq 2$ be integers.
	\begin{enumerate}
		\item \label{thm:1.1}
		An $m$-dimensional closed and connected manifold $M$ diffeomorphic to the total space of a smooth bundle over $S^n$ whose fiber is an {\rm (}$m-n${\rm )}-dimensional closed and connected manifold $F$ admits a round fold map $f:M \rightarrow {\mathbb{R}}^n$ having a globally trivial monodromy and for a suitable diffeomorphism $\phi:{\mathbb{R}}^n \rightarrow {\mathbb{R}}^n$ as in Definition \ref{def:2}, the preimage of the point ${\phi}^{-1}(0) \in {\mathbb{R}}^n$ is diffeomorphic to the disjoint union $F \sqcup F$ of two copies of the manifold $F$.
		\item \label{thm:1.2}
		 Let an $m$-dimensional closed and connected manifold $M$ admit a round fold map $f:M \rightarrow {\mathbb{R}}^n$.
		Let $M^{\prime}$ be the total space of an smooth bundle over $M$ whose fiber is a closed and connected manifold $F_1$. 
		Put the dimension of $M^{\prime}$ by $m^{\prime}$ and assume that $m^{\prime}>m$.
		Assume also that the following properties are enjoyed.
		\begin{enumerate}
			\item $f$ has a family of componentwisely trivial monodromies.
		 \item For $f$ and some small closed tubular neighborhood $N(C)$ of each connected component $C$ of $f(S(C))$, the restriction of the bundle $M^{\prime}$ over $M$ to $f^{-1}(N(C))$ is a trivial smooth bundle whose fiber is $F_1$.
		 \item For a suitable diffeomorphism ${\phi}_1:{\mathbb{R}}^n \rightarrow {\mathbb{R}}^n$ like $\phi$ in Definition \ref{def:2}, the preimage of the point ${{\phi}_1}^{-1}(0) \in {\mathbb{R}}^n$ is diffeomorphic to a closed smooth manifold $F_0$.
		 \end{enumerate}
	 Then $M^{\prime}$ admits a round fold map $f^{\prime}:M^{\prime} \rightarrow {\mathbb{R}}^n$ having a family of componentwisely trivial monodromies. Moreover, for a suitable diffeomorphism ${\phi}_2:{\mathbb{R}}^n \rightarrow {\mathbb{R}}^n$ like $\phi$ in Definition \ref{def:2}, the preimage of the point ${{\phi}_2}^{-1}(0) \in {\mathbb{R}}^n$ is diffeomorphic to the product $F_0 \times F_1$.
	\end{enumerate}
	\end{Thm}
Theorem \ref{thm:1} (\ref{thm:1.1}) is regarded as a specific case of (\ref{thm:1.2}). Consider a canonical projection of a unit sphere into the Euclidean space whose dimension is same.
\begin{Thm}[\cite{kitazawa0.1,kitazawa0.2,kitazawa0.5}]
	\label{thm:2}
	Let $m \geq n \geq 2$ be integers.
	An $m$-dimensional closed and connected manifold represented as a connected sum of the total spaces of smooth bundles over $S^n$ whose fibers are {\rm (}$m-n${\rm )}-dimensional standard spheres admits a round fold map $f$ into ${\mathbb{R}}^n$ enjoying the following properties where the connected sum is taken in the smooth category.
	\begin{enumerate}
		\item The index of each singular point is $0$ or $1$.
		\item The preimage of each regular value is the disjoint union of {\rm (}m-n{\rm )}-dimensional standard spheres.
		\item $f$ has a family of componentwisely trivial monodromies.
		
		\item The preimage of a regular value in the connected component of the regular value set of $f$ diffeomorphic to the interior of the  $n$-dimensional unit disk consists of $l+1$ connected components. By some suitable diffeomorphism $\phi:{\mathbb{R}}^n \rightarrow {\mathbb{R}}^n$ as in Definition \ref{def:2}, this connected component is regarded as the one containing the origin $0 \in {\mathbb{R}}^n$.  
	\end{enumerate} 
Furthermore, \cite{kitazawa0.2,kitazawa0.5} show that the converse is also true. 
\end{Thm}
For example, in the case $(m,n)=(4,2),(5,2)$, we have complete characterizations of closed and simply-connected manifolds admitting round fold maps considered in Theorem \ref{thm:2}.
We have stronger results in these cases. For example, we do not need to pose the property that the round fold map has a family of componentwisely trivial monodromies. For this, we need arguments on \cite{saekisuzuoka}, a kind of explicit theory on fold maps the preimages of regular values are disjoint unions of homotopy spheres, some well-known facts on diffeomorphism groups of some compact and connected manifolds (whose dimensions are at most $3$), and classifications of $5$-dimensional closed and simply-connected manifolds in \cite{barden}, for example. 

 Related to this, we present several known results on characterizations of manifolds admitting fold maps of fixed classes.

A {\it special generic} map means a fold map the index of whose singular point is always $0$. Morse functions with exactly two singular points on homotopy spheres are simplest examples and canonical projections of unit spheres are special generic. 

\begin{Thm}
	\label{thm:3}
	\begin{enumerate}
		\item \label{thm:3.1} {\rm (}Reeb's theorem{\rm )}
A closed and connected manifold $M$ admits a special generic map into $\mathbb{R}$ if and only if $M$ is either of the following two.
		\begin{enumerate}
		\item $M$ is a homotopy sphere of dimension $m \neq 4$.
		\item $M$ is a $4$-dimensional standard sphere. 
		\end{enumerate} 

\item
\label{thm:3.2}
 {\rm (}\cite{saeki}{\rm )} 
A closed and simply-connected manifold $M$ admits a special generic map into ${\mathbb{R}}^2$ if and only if $M$ is a manifold whose dimension is greater than $1$ in the previous case.
\item
\label{thm:3.3}
 {\rm (}\cite{calabi,saeki,saeki2}{\rm )}
A homotopy sphere of dimension $m>3$ admitting a special generic map into ${\mathbb{R}}^{m-j}$ must be a standard sphere for $j=1,2,3$.
\item
\label{thm:3.4}
 {\rm (}\cite{saeki}{\rm )} 
If An $m$-dimensional closed and simply-connected manifold $M$ admits a special generic map into the ${\mathbb{R}}^3$ for $m>3$, then $M$ is represented as a connected sum of the total spaces of smooth bundles over $S^2$
whose fibers are either of the following two where the connected sum is chosen in the smooth category.
\begin{enumerate}
	\item A homotopy sphere whose dimension is not $4$. 
	\item A $4$-dimensional standard sphere.
\end{enumerate}
Furthermore, in the case $m=4,5$ for example, the converse is also true.
\item
\label{thm:3.5}
 {\rm (}\cite{nishioka}{\rm )}
A $5$-dimensional closed and simply-connected manifold $M$ admits a special generic map into the ${\mathbb{R}}^4$ if and only if $M$ is a $5$-dimensional standard sphere or one represented as a connected sum of the total spaces of smooth bundles over $S^2$
whose fiber is $3$-dimensional standard spheres where the connected sum is chosen in the smooth category. Moreover, smooth bundles whose fibers are $3$-dimensional standard spheres here can be replaced by $3$-dimensional unit spheres.
	\end{enumerate}
	\end{Thm}

For example, we also have the following. They are closely related to Theorem \ref{thm:2} and Theorem \ref{thm:3} (\ref{thm:3.3}) for example.

\begin{Thm}
	\label{thm:4}
	Every $7$-dimensional oriented homotopy sphere $M$ admits a round fold map into ${\mathbb{R}}^4$ in Theorem \ref{thm:2}. Furthermore, we have a round fold map for the case $l \leq 2$.
	\begin{enumerate}
		\item $M$ admits such a round fold map satisfying $l=0$ if and only if it is a standard sphere. 
		\item $M$ admits such a round fold map satisfying $l=1$ if and only if it is a standard sphere or more generally, an oriented homotopy sphere of $16$ types of all $28$ types. In addition, standard spheres are of one of these $16$ types here. 
	\end{enumerate}
\end{Thm}
This is also due to structures of $7$-dimensional homotopy spheres. \cite{milnor} is a pioneering article, followed by \cite{eellskuiper}, for example.

A {\it graph manifold} is, in short, a $3$-dimensional closed and orientable manifold obtained by gluing the total spaces of smooth bundles over closed and connected surfaces which are orientable as $3$-dimensional manifolds along the boundaries by diffeomorphisms. 
\begin{Thm}[\cite{kitazawasaeki}]
	\label{thm:5}
	A $3$-dimensional closed and orientable manifold admits a round fold map into ${\mathbb{R}}^2$ if and only if it is a graph manifold.
	Furthermore, we can replace "round fold map" by "a round fold map having a family of componentwisely trivial monodromies".
	
	\end{Thm}
Note that most round fold maps
here have no globally trivial monodromies. See also \cite{kitazawasaeki2}. See also \cite{kitazawa0.6}, which will be revised drastically. Note also that this theorem is a stronger version of \cite{saeki3}.

\begin{Rem}
	\label{rem:2}
As another result of \cite{kitazawasaeki}, it has been announced that the converse of Theorem \ref{thm:2}, presented in the end there, is not true in the case $(m,n)=(3,2)$. Most of so-called {\it Lens spaces} and most of so-called {\it Seifert manifolds} are shown to be important examples.
\end{Rem}

Hereafer, we mainly consider {\it spin} manifolds. A {\it spin} manifold $X$ is an orientable smooth manifold the {\it $2$nd Stiefel-Whitney class} of whose tangent bundle is the zero element of $H^2(X;\mathbb{Z}/2\mathbb{Z})$, the 2nd homology group of $X$ whose {\it coefficient ring} is $\mathbb{Z}/2\mathbb{Z}$, the field of order $2$.

For any smooth manifold $X$, this is defined uniquely as an element of $H^2(X;\mathbb{Z}/2\mathbb{Z})$, the {\it 2nd homology group} of $X$ whose {\it coefficient ring} is $\mathbb{Z}/2\mathbb{Z}$, the field of order $2$. We explain about related notions shortly later again.

Homotopy spheres, the product of finitely many spin manifolds, and any manifold represented as a connected sum of finitely many spin manifolds where the connected sum is chosen in the smooth category are spin, for example.

The product of finitely many smooth manifolds at least one of which is not spin, and any manifold represented as a connected sum of finitely many manifolds at least one of which is not spin where the connected sum is chosen in the smooth category are not spin, for example. A closed smooth manifold which is homotopy equivalent to a closed spin manifold is also spin and one which is homotopy equivalent to a closed manifold which is not spin is not spin. 

The $k$-dimensional complex projective space ${\mathbb{C}P}^k$, which is a $2k$-dimensional smooth manifold, is spin if and only if $k$ is odd.

For more precise understanding on related notions, see \cite{milnorstasheff} again for example.
\begin{Thm}[\cite{kitazawa3}]
	\label{thm:6}
	There exists a family $\{M_{j}\}_{j \in \mathbb{Z}}$ of countably many $7$-dimensional closed and simply-connected oriented smooth manifolds
	enjoying the following properties.
	\begin{enumerate}
		\item The integral cohomology ring of $M_j$ is isomorphic to that of ${\mathbb{C}P}^2 \times S^3$.
		\item $M_j$ is spin.
		\item $M_{j}$ and $M_{-j}$ are homeomorphic but there exist no homeomorphisms preserving the orientations.
		\item $M_{j_1}$ and $M_{j_2}$ are not homeomorphic unless $j_1=\pm j_2$.
		\item $M_j$ admits a round fold map $f_j:M_j \rightarrow {\mathbb{R}}^4$ enjoying the following properties.
		\begin{enumerate}
			\item $f_j$ has a globally trivial monodromy.
			\item We set $f:=f_j$, $M:=M_j$ and $l=3$ and abuse the notation in Definition \ref{def:2}. Then for the preimages the following three hold.
			\begin{enumerate}
				\item The preimage of each point in ${\phi}^{-1}(\{x \in {\mathbb{R}}^n \mid 2<||x||<3\})$ is diffeomorphic to $S^3$.
				\item The preimage of each point in ${\phi}^{-1}(\{x \in {\mathbb{R}}^n \mid 1<||x||<2\})$ is diffeomorphic to $S^2 \times S^1$.
				\item The preimage of each point in ${\phi}^{-1}(\{x \in {\mathbb{R}}^n \mid 0 \leq |x||<1\})$ is diffeomorphic to the disjoint union $S^3 \sqcup (S^2 \times S^1)$ of $S^2 \times S^1$ and $S^3$.
			\end{enumerate}
		\end{enumerate}
	\end{enumerate}
	
	Furthermore, according to \cite{wang} with several methods of construction of round fold maps in \cite{kitazawa0.2,kitazawa0.5}, an arbitrary $7$-dimensional closed and simply-connected smooth manifold $M$ whose integral coefficient ring is isomorphic to that of ${\mathbb{C}P}^2 \times S^3$ and which is also spin admits a round fold map into ${\mathbb{R}}^4$. Note that $7$-dimensional manifolds here are not homeomorphic to ${\mathbb{C}P}^2 \times S^3$ since ${\mathbb{C}P}^2$ is not spin. 
\end{Thm}

The manifolds here and a complete classification of them are due to \cite{wang}. This preprint says that such $7$-dimensional manifolds are the total spaces of smooth bundles over ${\mathbb{C}P}^2$ whose fibers are diffeomorphic to $S^3$ and which are not trivial.

Our Main Theorem \ref{mthm:1} and related new works are closely related to this and this and its proof play important roles. 

\begin{Rem}
	\label{rem:3}
	Higher dimensional closed and simply-connected manifolds are ones whose dimensions are at least $5$. They are central objects in classical algebraic topology and differential topology. As presented in our introduction, they have been classified via sophisticated algebraic and abstract tools. $5$-dimensional ones have been classified by \cite{barden} for example. If the $2$nd integral homology groups of them are free, then they must be as ones in Theorems \ref{thm:2} and \ref{thm:3}. $6$-dimensional ones are classified through studies such as \cite{jupp,wall,zhubr} for example. It is somewhat difficult if there exists an ordered pair of elements of 2nd cohomology groups whose cup product is not the zero element of the $4$-th cohomology group. This is due to the relation $6=2 \times 3$ in short. For $7$-dimensional ones, more concrete classifications via concrete algebraic topology have been studied by \cite{crowleyescher,crowleynordstrom}, followed by \cite{kreck,wang}, for example.
\end{Rem}
\begin{Rem}
	\label{rem:4}
	Fold maps on manifolds in certain families of $7$-dimensional closed and simply-connected manifolds  whose integral cohomology rings are not isomorphic to that of ${\mathbb{C}P}^2 \times S^3$ are studied in \cite{kitazawa1,kitazawa2,kitazawa4}, for example. \cite{kitazawa1} also studies fold maps on such manifolds which may not be spin and manifolds whose dimensions are general. \cite{kitazawa1,kitazawa2} also concentrate on fold maps on closed (and simply-connected) manifolds of some classes which naturally extend the class of manifolds represented as connected sums of the total spaces of smooth bundles over standard spheres whose fibers are also standard spheres.
\end{Rem}
%Our Main Theorem \ref{mthm:2} may be regarded as a $6$-dimensional variant of this. Related to this, we give the definition of a ({\it generalized}) {\it Bott manifold}.

%We omit rigorous expositions on {\it Whitney sums} of vector bundles. See \cite{milnorstasheff} again. We also need {\it projectivizations} of vector bundles, which are obtained by considering bundles whose fibers are the projective spaces obtained canonically from the original fibers and whose structure groups are natural ones. We omit the exposition and Fundamental properties of the projectivizations of complex vector bundles are in \cite{choimasudasuh} for example.

%\begin{Def}
%	\label{def:4}
%	A smooth manifold $M$ is called a {\it generalized Bott manifold} if it is obtained in the following way.
%	\begin{enumerate}
%	\item Choose an arbitrary positive integer $k>0$ and a copy $M_0$ of the complex projective space ${\mathbb{C}}P^k$.
%	\item Consider a Whitney sum of finitely many complex vector bundles whose fibers are $1$-dimensional complex vector spaces. 
%	\item Take the total space of the projectivization of the bundle before and set $M_1$.
%	\item Consider a finite iteration of similar procedures one after another. In the $j$-th step $M_j$ denotes the resulting manifold where $j \geq 1$ is an integer. We choose $M:=M_{j_0}$ for some integer $j_0 \geq 1$.
%	\end{enumerate}
%Furthermore, if the Whitney sums are always chosen as ones of two bundles, then $M$ is called a {\it Bott manifold}.
%\end{Def}

	\section{Main Theorems.}
	Although systematic expositions are left to \cite{milnorstasheff} for example, we need the notions of the {\it $j$-th Stiefel-Whitney class} of a smooth manifold $M$ and the {\it $j$-th Pontrjagin class} of it and explain about them shortly.
	
	For any linear bundle over a base space $X$, we can define the {\it $j$-th Stiefel-Whitney class} as a uniquely defined element of $H^j(X;\mathbb{Z}/2\mathbb{Z})$. For any linear bundle over a base space $X$, we can define the {\it $j$-th Pontrjagin class} as a uniquely defined element of $H^j(X;\mathbb{Z})$.
	The {\it $j$-th Stiefel-Whitney class} of a smooth manifold $X$ is the uniquely defined element of $H^j(X;\mathbb{Z}/2\mathbb{Z})$, which is defined as the {\it $j$-th Stiefel-Whitney class} of the tangent bundle of $X$. Note that this is presented for the case $j=2$ before. This is a homotopical invariant for smooth manifolds.
	The {\it $j$-th Pontrjagin class} of $X$ is the uniquely defined element of $H^j(X;\mathbb{Z})$, which is defined as the {\it $j$-th Pontrjagin class} of the tangent bundle of $X$. This is an invariant for smooth manifolds. 
	%As an important fact, if the orientation is reversed, then only the sign of the $j$-th Pontrjagin class changed.
			
We prove Main Theorem \ref{mthm:1}. Some important methods are due to ones in the proof of Theorem \ref{thm:6}.

	\begin{proof}[A proof of Main Theorem \ref{mthm:1}]
	We can define a family $\{M_{0,j,s}\}_{(j,s) \in \mathbb{Z} \times \{0,1\}}$ of $m$-dimensional closed and simply-connected oriented manifolds.
	
	Let $S^2 \tilde{\times} S^2$ denote the total space of a linear bundle over $S^2$ whose fiber is the $2$-dimensional unit sphere and which is not trivial. This is well-known to be not spin and its integral cohomology ring is not isomorphic to $S^2 \times S^2$. We can argue as follows. For more precise arguments, see \cite{milnorstasheff} to know fundamental properties and arguments and \cite{wang} to know more explicit arguments which are also closely related to ours.
 
 	We can have a $7$-dimensional closed and simply-connected manfold $M_{0,j,s}$ diffeomorphic to the total space of a linear bundle over $S^2 \tilde{\times} S^2$ whose fiber is the ($m-4$)-dimensional unit sphere. The integral cohomology ring of $M_{0,j,s}$ is that of $(S^2 \tilde{\times} S^2) \times S^3$.

  We can find an element $h_{\rm b} \in H_2(S^2 \tilde{\times} S^2;\mathbb{Z})$ of the base space represented by a $2$-dimensional standard sphere enjoying the following properties and arguments.
 
 . 
 \begin{itemize}
 	\item The sphere $S_{\rm b}$ is regarded as the image of the section of the bundle $S^2 \tilde{\times} S^2$ over $S^2$, whose fiber is the $2$-dimensional unit sphere.  
 	\item We can find the element $h_{\rm f} \in H_2(S^2 \tilde{\times} S^2;\mathbb{Z})$ represented by the fiber of the bundle over $S^2$ and $\{h_{\rm b},h_{\rm f}\}$ is a basis of $H_2(S^2 \tilde{\times} S^2;\mathbb{Z})$.
 	\item For the basis $\{h_{\rm b},h_{\rm f}\}$, we can define the dual ${h_{\rm b}}^{\ast} \in H^2(S^2 \tilde{\times} S^2;\mathbb{Z})$ to $h_{\rm b}$ and the dual ${h_{\rm f}}^{\ast}\in H^2(S^2 \tilde{\times} S^2;\mathbb{Z})$ to $h_{\rm f}$. Furthermore, the cup product ${h_{\rm f}}^{\ast} \cup {h_{\rm f}}^{\ast}$ is a generator of $H^4(S^2 \tilde{\times} S^2;\mathbb{Z})$, isomorphic to $\mathbb{Z}$, the cup product ${h_{\rm b}}^{\ast} \cup {h_{\rm f}}^{\ast}$ is a generator of $H^4(S^2 \tilde{\times} S^2;\mathbb{Z})$, and the cup product ${h_{\rm b}}^{\ast} \cup {h_{\rm b}}^{\ast}$ is the zero element of $H^4(S^2 \tilde{\times} S^2;\mathbb{Z})$. 
 	%\item The Poincar\'e dual ${\rm PD}({h_{\rm f}}^{\ast})$ to ${h_{\rm f}}^{\ast}$ is $h_{\rm b}$. The Poincar\'e dual ${\rm PD}({h_{\rm b}}^{\ast})$ to ${h_{\rm b}}^{\ast}$ is $h_{\rm b}-h_{\rm f}$.
 \end{itemize}	
We have a desired family of the manifolds with related arguments as follows. For detailed arguments on ($2$nd) Stiefel-Whitney classes and ($1$st) Pontrjagin classes, consult \cite{milnorstasheff} and \cite{wang} again.
\begin{itemize}
	\item We can have $M_{0,j,0}$ as a manifold whose $1$st Pontrjagin class is $4j$ times a generator of $H^4(M_{0,j,0};\mathbb{Z})$, isomorphic to $\mathbb{Z}$, and which is not spin.
	
	We consider the trivial linear bundle over $S^2 \tilde{\times} S^2$ whose fiber is the $3$-dimensional unit sphere $S^3$ and on a suitable small $4$-dimensional copy of the unit disk $D^4$ smoothly embedded here, we exchange the way of attaching the trivial linear bundle over the embedded copy of the $4$-dimensional unit disk $D^4$ whose fiber is the $3$-dimensional unit sphere $S^3$. This contributes to "$4j$" before. 

\item We can have $M_{0,j,1}$ as a manifold whose $1$st Pontrjagin class is $4j$ times a generator of $H^4(M_{0,j,0};\mathbb{Z})$, isomorphic to $\mathbb{Z}$, and which is spin.

We consider the trivial linear bundle over $S^2 \tilde{\times} S^2$ whose fiber is the $3$-dimensional unit sphere $S^3$ and on a suitable small moothly embedded manifold diffeomorphic to a copy of the $4$-dimensional unit disk $D^4$ smoothly embedded here, we exchange the way of attaching the trivial linear bundle over the $4$-dimensional manifold whose fiber is the $3$-dimensional unit sphere $S^3$. This contributes to "$4j$" before. 
We can choose a small closed tubular neighborhood of the $2$-dimensional sphere $S_{\rm b}$.
 We exchange the way of attaching the trivial linear bundle over the closed tubular neighborhood whose fiber is the $3$-dimensional unit sphere $S^3$ to change our $7$-dimensional closed and simply-connected manifold into one whose 2nd Stiefel-Whitney-class is the zero element. 

This contributes to "$1$" in "$4j+1$" before.

\end{itemize}
We can regard each manifold $M_{0,j,s}$ as the total space of a smooth bundle over $S^2$ in a natural way. If $s=0$ ($s=1$), then the fiber is regarded as a manifold diffeomorphic to the total space of a trivial linear bundle over $S^2$ whose fiber is the ($m-4$)-dimensional unit sphere $S^{m-4}$. 
%Note that the total space of a linear bundle over $S^2$ whose fiber is the ($m-4$)-dimensional unit sphere $S^{m-4}$ is well-known to be (not) spin if it is (resp. not) trivial.

We can exchange each projection into a fold map on another $7$-dimensional closed and simply-connected manifold $M_{1,j,s}$ onto $S^2$ whose singular set is a circle and the restriction of which to the singular set is an embedding. Furthermore, we can do this in such a way that the preimage of a regular value in one of the connected components of the regular value set is, an ($m-2$)-dimensional standard sphere and the preimage of a regular value in the other one of the connected components of the regular value set is, diffeomorphic to the original total space of the originally considered trivial linear bundle over $S^2$ whose fiber is the ($m-4$)-dimensional unit sphere $S^{m-4}$. 

We present construction of a round fold map on $M_{1,j,s}$.
We decompose the manifold into two $m$-dimensional compact manifolds.
\begin{itemize}
	\item The preimage of the disjoint union of two copies of the $4$-dimensional unit disk $D^4$ smoothly embedded into distinct connected components of the regular value set. Let $\tilde{M_{1,j,s}}_1$ denote this manifold. 
	\item The complementary set of the interior of the previous manifold. Let $\tilde{M_{1,j,s}}_2$ denote this manifold.
\end{itemize}

In a natural way, $\tilde{M_{1,j,s}}_1$ is regarded as the total space of a trivial smooth bundle over ${D^2}_r:={x \in \{\mathbb{R}}^2 \mid ||x|| \leq r\}$ for $r>0$ whose fiber is diffeomorphic to the disjoint union of the total space of the bundle over $S^2$ whose fiber is the ($m-4$)-dimensional unit sphere $S^{m-4}$ and a copy of the ($m-2$)-dimensional unit sphere $S^{m-2}$. $\tilde{M_{1,j,s}}_2$ admits a map regarded as the product map of a Morse function enjoying the following properties and the identity map on the boundary $\partial {D^2}_r \subset {D^2}_r$.
\begin{itemize}
	\item The manifold of the domain of the Morse function is a manifold obtained by removing the interior of a copy of the ($m-1$)-dimensional unit disk $D^{m-1}$ smoothly embedded in the interior of the total space of a trivial linear bundle over $S^2$ whose fiber is the ($m-3$)-dimensional unit disk $D^{m-3}$.
	\item The singular points are in the interior of the manifold of the domain and at distinct singular points, the values are distinct.
	\item The function has exactly three singular points.
	\end{itemize}

We can glue the projection of the trivial bundle on $\tilde{M_{1,j,s}}_1$ and the product map on $\tilde{M_{1,j,s}}_2$ suitable way to obtain a round fold map from $M_{1,j,s}$ into ${\mathbb{R}}^2$.

The singular set consists of three connected components. By considering the composition with a natural projection onto $\mathbb{R}$ in a suitable way, we have a Morse function with exactly six singular points.

In addition, we can do our construction of new fold maps in such a way that we have some important topological properties of the resulting $m$-dimensional manifold $M_{1,j,s}$ as follows.
	\begin{itemize}
		\item There exists an element of $e \in H_2(M_{1,j,s};\mathbb{Z})$ which is not divided by any integer greater than $1$ and which is a generator of the free group $H_2(M_{1,j,s};\mathbb{Z})$.
		\item We can define the dual $e^{\ast}$ to this due to the previous fact that $H_2(M_{1,j,s};\mathbb{Z})$ is free and of rank $1$.
		\item The Poincar\'e dual ${\rm PD}(e^{\ast})$ to the dual $e^{\ast}$, which is an element of $H^2(M_{1,j,s};\mathbb{Z})$, is a generator of $H_{m-2}(M_{1,j,s};\mathbb{Z})$, and represented by the preimage of each regular value.
		\item The Poincar\'e dual to the square $e^{\ast} \cup e^{\ast} \in H^4(M_{1,j,s};\mathbb{Z})$ is represented by a suitable submanifold of the preimage of a regular value. For the suitable submanifold, for example, we can choose one which is originally a fiber of the total space of the original smooth bundle over $S^2 \tilde{\times} S^2$. Remember that this fiber is diffeomorphic to $S^2 \times S^{m-4}$. Let $e^{\prime} \in H_{m-4}(M_{1,j,s};\mathbb{Z})$ denote the Poincar\'e dual. $e^{\ast} \cup e^{\ast} \in H^4(M_{1,j,s};\mathbb{Z}) \cong \mathbb{Z}$ is also a generator of the $4$th integral homology group.
        \item The ($m-4$)-th integral homology group $H_{m-4}(M_{1,j,s};\mathbb{Z})$ is free and of rank $1$ if $m-4 \neq 4$ and rank $2$ if $m-4=4$. 
        We can define the dual to the element $e^{\prime} \in H_{m-4}(M_{1,j,s};\mathbb{Z})$, defined before, by taking a natural basis of $H_{m-4}(M_{1,j,s};\mathbb{Z})$ (both in the cases $m-4 \neq 4$ and $m-4=4$).
        The cup product of the square $e^{\ast} \cup e^{\ast} \in H^4(M_{1,j,s};\mathbb{Z})$ and the dual ${e^{\prime}}^{\ast} \in H^{m-4}(M_{1,j,s};\mathbb{Z})$ to the element $e^{\prime} \in H_{m-4}(M_{1,j,s};\mathbb{Z})$, which is defined before, is a generator of $H^m(M_{1,j,s};\mathbb{Z})$. 
         \item The following cup products are the zero elements. We can know them due to Poincare duality, or more precisely, intersection theory for closed submanifolds with no boundaries and elements of integral homology groups represented by these submanifolds.  
         \begin{itemize}
         	\item The cup product of an element of $H^2(M_{1,j,s};\mathbb{Z})$ and an element of $H^4(M_{1,j,s};\mathbb{Z})$. 
         	\item The cup product of an element of $H^2(M_{1,j,s};\mathbb{Z})$ and an element of $H^{m-4}(M_{1,j,s};\mathbb{Z})$.
         	\item The square $(e^{\ast} \cup e^{\ast}) \cup (e^{\ast} \cup e^{\ast}) \in H^8(M_{1,j,s};\mathbb{Z})$ of the square $e^{\ast} \cup e^{\ast} \in H^4(M_{1,j,s};\mathbb{Z})$ is the zero element.
         	\item The square ${e^{\prime}}^{\ast} \cup {e^{\prime}}^{\ast} \in H^8(M_{1,j,s};\mathbb{Z})$ is the zero element.
         	\end{itemize}
         \item For the 2nd Stiefel-Whitney classes and the 1st Pontrjagin classes of our new $7$-dimensional closed and simply-connected manifolds, the following two hold.
         \begin{itemize}
         	
         	\item We can have $M_{1,j,0}$ as a manifold whose $1$st Pontrjagin class is $4j-1$ times a generator $e^{\ast} \cup e^{\ast} \in H^4(M_{1,j,0};\mathbb{Z})$ and which is not spin.

         	\item We can have $M_{1,j,1}$ as a manifold whose $1$st Pontrjagin class is $4j$ times a generator $e^{\ast} \cup e^{\ast} \in H^4(M_{1,j,1};\mathbb{Z})$ and which is spin.

         \end{itemize}
         
         For related arguments on Pontrjagin classes, consult \cite{milnorstasheff} and \cite{wang} again.
\end{itemize}

This argument together with the number of singular points of the Morse function enables us to see that the integral cohomology ring of $M_{1,j,s}$ is isomorphic to that of ${\mathbb{C}P}^2 \times S^{m-4}$ and that the manifold is simply-connected.

Now we have a desired round fold map on a desired manifold $M_{j,s}:=M_{1,j,s}$.   
  	\end{proof}
  In the following, Theorem \ref{thm:7} summarizes our construction in the proof before.
  \begin{Thm}
  	\label{thm:7}
  	In Main Theorem \ref{mthm:1}, we can choose the manifolds as follows.
  	\begin{enumerate}
  		\item The 1st Pontrjagin class of $M_{j,s}$ is $4j+s-1$ times a generator of $H^4(M_{j,s};\mathbb{Z}) \cong \mathbb{Z}$. 
  		\item $M_{j,s}$ is {\rm (}not{\rm )} spin if $s=1$ {\rm (}resp. $s=0${\rm )}.
  	\end{enumerate}
  Furthermore, round fold maps are constructed as ones having globally trivial monodromies.
  \end{Thm}

  This is regarded as a new theorem closely related to Theorem \ref{thm:6} in several senses. For example, we only had round fold maps into ${\mathbb{R}}^4$ for the case where the manifolds are spin and $7$-dimensional.

 \begin{Rem}
	\label{rem:5}
	Let $m \geq 8$ be an even integer.
	 The Euler number of an arbitrary $m$-dimensional manifold in Main Theorem \ref{mthm:1} is $2$.
	 
	  Suppose that the manifold has the structure of some smooth bundle over $S^2$ and $m$ is divisible by $4$. In such a case, the Euler number of $S^2$ is $2$ and the Euler number of an arbitrary ($m-2$)-dimensional closed orientable manifold must be even. These facts contradict the fact that the Euler number of the given $m$-dimensional manifold is $2$ and not divisible by $4$.
	  
	   Suppose that the manifold has the structure of some smooth bundle over $S^2$ and $m$ is not divisible by $4$. In such a case, the Euler number of $S^2$ is $2$ and the Euler number of an ($m-2$)-dimensional closed orientable manifold $M^{\prime}$ may be odd. There exists an element $u \in H^{\frac{m-2}{2}}(M^{\prime};\mathbb{Z})$ such that the square $u \cup u \in H^{m-2}(M^{\prime};\mathbb{Z})$ is not divisible by $2$. $\frac{m-2}{2} \geq 3$ holds. This contradicts the structure of the integral cohomology ring of the given $m$-dimensional manifold.

\end{Rem}
In the following, Theorem \ref{thm:8} summarizes Remark \ref{rem:5}.
\begin{Thm}
	\label{thm:8}
	Let $m \geq 8$ be an even integer.	
	We can not obtain round fold maps on the $m$-dimensional manifolds in Main Theorem \ref{mthm:1} directly via Theorem \ref{thm:1} {\rm (}\ref{thm:1.1}{\rm )}.
\end{Thm}
 \begin{Rem}
	\label{rem:6}
We have some additional arguments closely related to Remark \ref{rem:5}.
Let $m=7$ in Main Theorem \ref{mthm:1}. 

We can see that these manifolds do not have the structures of smooth bundles over $S^4$, $S^5$ or $S^6$. In fact, by a fundamental argument on homotopy groups, if the manifold has this structure over $S^4$ or $S^5$, then the fiber must be simply-connected and as a result diffeomorphic to a sphere. We can easily see that this is a contradiction. In the case where the base space is $S^6$, we can argue similarly. We cannot construct round fold maps there into ${\mathbb{R}}^n$ for $n=4,5,6$ via Theorem \ref{thm:1} (\ref{thm:1.1}).  

If the manifold of dimension $m=7$ is spin or the Pontrjagin class of the manifold is not $3 $ times a generator of the 4th integral cohomology group, then it does not have the structure of a smooth bundle over $S^3$. We explain about this.

Suppose that this has the structure of a smooth bundle over $S^3$. Then the fiber is a $4$-dimensional closed and simply-connected manifold, which follows from a fundamental argument on homotopy groups. By an elementary argument on (co)homology groups, we can see that the 2nd integral homology group of the fiber is free and of rank $1$. This is homeomorphic to the $2$-dimensional projective space ${\mathbb{C}P}^2$ by virtue of important theory of $4$-dimensional manifolds, discussed in \cite{freedman} for example. This is not spin, by virtue of the theory. The 1st Pontrjagin class of ${\mathbb{C}P}^2$ is $3$ times a generator of its $4$th integral cohomology group. This is a kind of important facts on ($1$st) Pontrjagin classes, discussed in \cite{milnorstasheff} for example. 
This is a contradiction. 

We can see that we cannot construct round fold maps there into ${\mathbb{R}}^3$ via Theorem \ref{thm:1} (\ref{thm:1.1}). Note that the product $S^3 \times {\mathbb{C}P}^2$ admits a round fold map into ${\mathbb{R}}^4$. In fact, we can consider a Morse function on ${\mathbb{C}P}^2$ and the product map of this and the identity map on $S^3$ and we can smoothly embed the target of the resulting map into ${\mathbb{R}}^4$ to have a round fold map into ${\mathbb{R}}^4$. The fact on the existence of round fold maps in Remark \ref{rem:1} gives round fold maps into ${\mathbb{R}}^n$ for $n={\rm (}1,{\rm )} 2, 3$ there.
\end{Rem}
 \begin{Rem}
	\label{rem:7}
Related to Remarks \ref{rem:5} and \ref{rem:6}, we do not know whether we can construct round fold maps on the manifolds in Main Theorem \ref{mthm:1} via Theorem \ref{thm:1} (\ref{thm:1.2}) in general. Remember that Theorem \ref{thm:1} (\ref{thm:1.2}) is extends Theorem \ref{thm:1} (\ref{thm:1.1}) for example.
\end{Rem}

We close this section by some general remarks on round fold maps.

\begin{Rem}
	\label{rem:8}
	A round fold map into ${\mathbb{R}}^2$ always presents the structure of an {\it open book} of the manifold of the domain. It is important that in considerable cases closed and simply-connected manifolds have such structures and in most of such cases we can also have nicer ones according to \cite{winkelnkemper}.
	For example, if a closed and simply-connected manifold of the dimension is not divisible by $4$, then it has the structure of a nice open book from the viewpoint of the (integral) homology groups. It is an open problem whether for an arbitrary open book, we can have a round fold map presenting the desired open book in the canonical way. 
	
	For example, some theory of $3$-dimensional manifolds presents interesting examples.
	So-called {\it hyperbolic} manifolds form a wide class of $3$-dimensional manifolds and manifolds whose dimensions are general. Graph manifolds are never hyperbolic, for example.
	A $3$-dimensional closed {\it hyperbolic} manifold has the structures of open books. However, it does not have round fold maps, or fold maps such that the restrictions to the singular sets are embeddings, due to \cite{saeki3}.   
\end{Rem}
\begin{Rem}
	\label{rem:9}
	{\it Stable} maps form an important class of smooth maps in singularity theory of smooth maps and applications to geometry of manifolds. For example, Morse functions whose singular values are always distinct at distinct singular points and fold maps  whose restrictions to the singular sets are embeddings are stable.

	Stable maps	exist densely if the pair of the dimensions of the smooth (closed) manifold of the domain and the smooth manifold with no boundary of the target is nice. As the topology of the set of (all) smooth maps, we consider the so-called {\it $C^{\infty}$ Whitney topology}. For example, if the dimension of the manifold of the target is smaller than $6$, then the pair of the dimensions is nice. 

	For example, it is well-known that
	a smooth function on a closed manifold is stable if and only it is a Morse function as before. According to \cite{thom,whitney}, the singular set of a stable map into ${\mathbb{R}}^2$ on a closed manifold whose dimension is greater than or equal to $2$ is a $1$-dimensional smooth closed and compact submanifold with no boundary and the restriction there is a topological (PL or piecewwise smooth) immersion. Furthermore, the restriction of the map to the complementary set of finitely many points in the singular set is a fold map. In addition, these finitely many points are singular points called {\it cusp} points. A cusp point is defined as a singular point where the smooth map has forms of a certain type.
	It has been also shown that a smooth map on a closed manifold whose dimension is greater than or equal to $2$ into ${\mathbb{R}}^2$ is stable if the following properties are enjoyed.
	\begin{enumerate}
				\item The singular set is a $1$-dimensional smooth closed and compact submanifold with no boundary.
				\item Except finitely many singular points, around each singular point, it is a fold map. These finitely many singular points are cusp points.
				\item The restriction to the singular set is a topological (PL or piecewise smooth) immersion. Furthermore, the preimage of the immersion has at most two points. If the preimage of the point in the space of the target consists of exactly two points, then the points of the preimage are not cusp points and the sum of the images of the differentials there and the tangent vector space at the point in the space of the target agree.
				\end{enumerate}

		For systematic theory, see \cite{golubitskyguillemin}. 
	
	Recently, for example, \cite{baykursaeki,baykursaeki2} show that any $4$-dimensional closed and orientable smooth manifold admits a stable map into ${\mathbb{R}}^2$ enjoying good differential topological properties. More precisely, the restriction to the singular set, which is a $1$-dimensional smooth closed and compact submanifold with no boundary, is a topological (PL or piecewwise smooth) embedding.
\end{Rem}
\begin{Rem}
	\label{rem:10}
Recently, in the conference \\

 https://www.fit.ac.jp/~fukunaga/conf/sing202206.html, \\
 
\noindent Saeki has announced a result which is regarded as a higher dimensional variant of the result for stable maps on $4$-dimensional closed manifolds into ${\mathbb{R}}^2$ in Remark \ref{rem:5}. More precisely, he has announced a variant for stable maps of closed manifolds whose dimensions are arbitrary and greater than $3$ into ${\mathbb{R}}^2$.  

\end{Rem}

\section{Declarations.}
\label{sec:4}
\thanks{The author is a member of JSPS KAKENHI Grant Number JP17H06128 "Innovative research of geometric topology and singularities of differentiable mappings" (Principal Investigator: Osamu Saeki).}

%The author would like to thank anonymous referees for important and interesting comments on submitted versions of the paper. This has improved the paper and the author.

We declare that we have no associated data on the present paper. 


\begin{thebibliography}{30}
\bibitem{kitazawa4} N. Kitazawa, \textsl{$7$-dimensional closed simply-connected and spin manifolds having 2nd integral cohomology classes whose squares are not divisible by $2$ and stable fold maps on them}, arxiv:2104.10871.
\end{thebibliography}
\end{document}